\documentclass[10pt, a4paper]{article}

\usepackage{amsfonts, latexsym, amsmath, amssymb, fancyhdr, amsthm}
\usepackage[english]{babel}

\newtheorem{theorem}{Theorem}[section]
\newtheorem{lemma}{Lemma}[section]
\newtheorem{proposition}{Proposition}[section]
\newtheorem{definition}{Definition}[section]
\newtheorem{corollary}{Corollary}[section]
\newtheorem{remark}{Remark}

\title{\textbf{On the stability of a nonlinear maturity structured model of cellular proliferation\thanks{This paper has been published in Dis. Cont. Dyn. Sys. Ser. A, 12 (3), 501-522, 2005.}}}
\author{Mostafa Adimy$^\ddagger$,\quad Fabien Crauste$^\ddagger$ \quad and \quad Laurent Pujo-Menjouet$^\dagger$}
\date{Year 2004}

\begin{document}

\maketitle

\begin{center}
{\large $^\ddagger$} \emph{Laboratoire de Math\'ematiques Appliqu\'ees, FRE 2570,} \\ \emph{Universit\'e de Pau et des Pays de l'Adour, Avenue de l'universit\'e, 64000 Pau, France.}\\\emph{E-mail: mostafa.adimy@univ-pau.fr}, \quad \emph{E-mail: fabien.crauste@univ-pau.fr}\\
\quad\\
{\large $^\dagger$} \emph{Department of Physiology, McGill
University, McIntyre Medical Sciences Building,}\\ \emph{3655
Promenade Sir William Osler, Montreal, QC, Canada H3G
1Y6.}\\\emph{ E-mail: pujo@cnd.mcgill.ca}
\end{center}

\quad

\begin{abstract}
We analyse the asymptotic behaviour of a nonlinear mathematical
model of cellular proliferation which describes the production of
blood cells in the bone marrow. This model takes the form of a
system of two maturity structured partial differential equations,
with a retardation of the maturation variable and a time delay
depending on this maturity. We show that the stability of this
system depends strongly on the behaviour of the immature cells
population. We obtain conditions for the global stability and the
instability of the trivial solution.
\end{abstract}

\bigskip{}


\section{Introduction and motivation}

This paper is devoted to the analysis of a maturity structured
model which involves descriptions of process of blood production
in the bone marrow (hematopoiesis). Cell biologists recognize two
main stages in the process of hematopoietic cells: a resting stage
and a proliferating stage (see Burns and Tannock \cite{burns}).

The resting phase, or $G_0$-phase, is a quiescent stage in the
cellular development. Resting cells mature but they can not
divide. They can enter the proliferating phase, provided that they
do not die. The proliferating phase is the active part of the
cellular development. As soon as cells enter the proliferating
phase, they are committed to divide, during mitosis. After
division, each cell gives birth to two daughter cells which enter
immediatly the resting phase, and complete the cycle.
Proliferating cells can also die without ending the cycle.

The model considered in this paper has been previously studied by
Mackey and Rudnicki in 1994 \cite{mackey1994} and in 1999
\cite{mackey1999}, in the particular case when the proliferating
phase duration is constant. That is, when it is supposed that all
cells divide exactly at the same age. Numerically, Mackey and Rey
\cite{mackey1995_1,mackey1995_2}, in 1995, and Crabb \emph{et al.}
\cite{crabb1996_1,crabb1996_2}, in 1996, obtained similar results
as in \cite{mackey1994}. The model in \cite{mackey1994} has also
been studied by Dyson \emph{et al} \cite{webb1996} in 1996 and
Adimy and Pujo-Menjouet \cite{adimypujo,adimypujo2} in 2001 and
2003, but only in the above-mentioned case. These authors showed
that the uniqueness of the entire population depends, for a finite
time, only on the population of small maturity cells.

However, it is believed that, in the most general situation in
hematopoiesis, all cells do not divide at the same age (see
Bradford \emph{et al.} \cite{bradford}). For example, pluripotent
stem cells (the less mature cells) divide faster than committed
stem cells (the more mature cells).

Mackey and Rey \cite{mackey1993}, in 1993, considered a model in
which the time required for a cell to divide is not identical
between cells, and, in fact, is distributed according to a
density. However, the authors made only a numerical analysis of
their model. Dyson \emph{et al.} \cite{webb2000, webb2000_2}, in
2000, also considered an equation in which all cells do not divide
at the same age. But they considered only one phase (the
proliferating one) which does not take into account the
intermediary flux between the two phases. Adimy and Crauste
\cite{adimycrauste2003}, in 2003, studied a model in which the
proliferating phase duration is distributed according to a density
with compact support. The authors proved local and global
stability results.

In \cite{adimycrauste}, Adimy and Crauste developed a mathematical
model of hematopoietic cells population in which the time spent by
each cell in the proliferating phase, before mitosis, depends on
its maturity at the point of commitment. More exactly, a cell
entering the proliferating phase with a maturity $m$ is supposed
to divide a time $\tau=\tau(m)$ later. This hypothesis can be
found, for example, in Mitchison \cite{mitchison} (1971) and John
\cite{john} (1981), and, to our knowledge, it has never been used,
except by Adimy and Pujo-Menjouet in \cite{adimypujo2003}, where
the authors considered only a linear case. The model obtained in
\cite{adimycrauste} is a system of nonlinear first order partial
differential equations, with a time delay depending on the
maturity and a retardation of the maturation variable. The basic
theory of existence, uniqueness, positivity and local stability of
this model was investigated.

Many cell biologists assert that the behaviour of immature cells
population is an important consideration in the description of the
behaviour of full cells population. The purpose of the present
work is to analyse mathematically this phenomenon in our model. We
show that, under the assumption that cells, in the proliferating
phase, have enough time  to divide, that is, $\tau(m)$ is large
enough, then the uniqueness of the entire population depends
strongly, for a finite time, on the population with small
maturity. This result allows us, for example, to describe the
destruction of the cells population when the population of small
maturity cells is affected (see Corollary
\ref{coraplascticanemia}).

In \cite{mackey1999}, Mackey and Rudnicki provided a criterion for
global stability of their model. However, these authors considered
only the case when the mortality rates and the rate of returning
in the proliferating cycle are independent of the maturity
variable. Thus, their criterion can not be applied directly to our
situation.

This paper extends some local analysis of Adimy and Crauste
\cite{adimycrauste} to global results. It proves the connection
between the global behaviour of our model and the behaviour of
immature cells ($m=0$).

The paper is organised as follows. In the next section, we present
the equations of our model and we give an integrated formulation
of the problem, by using the semigroup theory. In section
\ref{sectionuniquenessresult}, we show an uniqueness result which
stresses the dependence of the entire population with small
maturity cells population. In Section
\ref{sectionbehaviourimmaturecells}, we focus on the behaviour of
the immature cells population, which satisfies a system of delay
differential equations. We study the stability of this system by
using a Lyapunov functionnal. In Section
\ref{asymptoticbehavioursection}, we prove that the global
stability of our model depends on its local stability and on the
stability of the immature cells population. Finally, in Section
\ref{sectioninstability}, we give an instability result.

\section{Equations of the model and integrated formulation}

Let $N(t,m)$ and $P(t,m)$ denote, respectively, the population
densities of resting and proliferating cells, at time $t$ and with
a maturity level $m$.

The maturity is a continuous variable which represents what
composes a cell, such as proteins or other elements one can
measure experimentally. It is supposed to range,  in the two
phases, from $m=0$ to $m=1$. Cells with maturity $m=0$ are the
most primitive stem cells, also called immature cells, whereas
cells with maturity $m=1$ are ready to enter the bloodstream, they
have reached the end of their development.

In the two phases, cells mature with a velocity $V(m)$, which is
assumed to be continuously differentiable on $\left[ 0,1\right] $,
positive on $\left( 0,1\right] $ and such that $V(0)=0$ and
\begin{equation}
\displaystyle\int_{0}^{m}\dfrac{ds}{V(s)}=+\infty , \qquad \text{
for }m\in \left( 0,1\right].  \label{(H.1).1}
\end{equation}
Since $\int_{m_{1}}^{m_{2}}\frac{ds}{V(s)}$, with $m_1<m_2$, is
the time required for a cell with maturity $m_{1}$ to reach the
maturity $m_{2}$, then Condition (\ref{(H.1).1}) means that a cell
with very small maturity needs a long time to become mature.

\noindent For example, Condition (\ref{(H.1).1}) is satisfied if
\begin{displaymath}
V(m)\underset{m\to 0}{\sim} \alpha m^p,\quad \textrm{ with }
\alpha>0 \textrm{ and } p\geq 1.
\end{displaymath}

In the resting phase, cells can die at a rate $\delta=\delta(m)$
and can also be introduced in the proliferating phase with a rate
$\beta$. In the proliferating phase, cells can also die, by
apoptosis (a programmed cell death), at a rate $\gamma=\gamma(m)$.
The functions $\delta$ and $\gamma$ are supposed to be continuous
and nonnegative on $[0,1]$. The rate $\beta$ of re-entry in the
proliferating phase is supposed to depend on cells maturity and on
the resting population density (see Sachs \cite{sachs}), that is,
$\beta=\beta(m,N(t,m))$. The mapping $\beta$ is supposed to be
continuous and positive.

Proliferating cells are committed to undergo mitosis a time $\tau$
after their entrance in this phase. We assume that $\tau$ depends
on the maturity of the cell when it enters the proliferating
phase, that means, if a cell enters the proliferating phase with a
maturity $m$, then it will divide a time $\tau=\tau(m)$ later.

\noindent The function $\tau$ is supposed to be positive,
continuous on $[0,1]$, continuously differentiable on $(0,1]$ and
such that
\begin{equation}\label{condtau}
\tau^{\prime }(m)+\frac{1}{V(m)}>0, \qquad \textrm{ for }
m\in(0,1].
\end{equation}
One can notice that this condition is always satisfied in a
neighborhood of the origin, because $V(0)=0$, and is satisfied if
we assume, for example, that $\tau$ is increasing (which describes
the fact that the less mature cells divide faster than more mature
cells).

\noindent Under Condition (\ref{condtau}), if $m\in(0,1]$ is
given, then the mapping
\begin{displaymath}
\widetilde{m} \mapsto
\int_{\widetilde{m}}^{m}\frac{ds}{V(s)}-\tau(\widetilde{m})
\end{displaymath}
is continuous and strictly decreasing from $(0,m]$ into
$[-\tau(m),+\infty)$. Hence, we can define a function $\Theta :
(0,1]\to (0,1]$, by
\begin{displaymath}
\int_{\Theta(m)}^m \frac{ds}{V(s)}=\tau(\Theta(m)), \quad \textrm{
for } m\in (0,1].
\end{displaymath}
The quantity $\Theta (m)$ represents the maturity of a cell at the
point of commitment when this cell divides at a maturity level
$m$. The function $\Theta$ is continuously differentiable and
strictly increasing on $(0,1]$ and satisfies
\begin{displaymath}
\lim_{m\to 0} \Theta(m)=0 \qquad \textrm{ and } \qquad 0<\Theta
(m)<m, \quad \textrm{ for } m\in(0,1].
\end{displaymath}

\noindent If we consider the characteristic curves $\chi :
(-\infty,0]\times[0,1] \to [0,1]$, solutions of the ordinary
differential equation
\begin{displaymath}
\left\{ \begin{array}{rcll}
\displaystyle\frac{d\chi}{ds}(s,m) & = & V(\chi(s,m)),& s \leq 0 \textrm{ and } m\in[0,1], \\
\chi(0,m)        & = & m, &
\end{array} \right.
\end{displaymath}
then, it is easy to check that, for $m\in[0,1]$, $\Theta (m)$ is
the unique solution of the equation
\begin{equation} \label{propertyoftheta}
x=\chi(-\tau(x),m).
\end{equation}
The characteristic curves $\chi(s,m)$ represent the evolution of
the cell maturity to reach a maturity $m$ at time $0$ from a time
$s \leq 0$. They satisfy $\chi(s,0)=0$ and $\chi(s,m) \in (0,1]$
for $s\leq 0$ and $m \in (0,1]$. Moreover, we can verify that the
characteristic curves are given by
\begin{equation}\label{charach}
\chi(s,m)=h^{-1}(h(m)e^{s}), \quad \textrm{ for } s\leq 0 \textrm{
and } m\in[0,1],
\end{equation}
where the continuous function $h :[0,1] \to [0,1]$ is defined by
\begin{displaymath}
h(m)=\left\{ \begin{array}{ll}
\exp\bigg(-\displaystyle\int_{m}^{1} \frac{ds}{V(s)} \bigg),&\textrm{for } m \in (0,1],\\
0,& \textrm{for } m=0.
\end{array} \right.
\end{displaymath}
Since $h$ is increasing, the two functions $s\mapsto \chi(s,m)$
and $m\mapsto \chi(s,m)$ are also increasing.

At the end of the proliferating phase, a cell with a maturity $m$
divides into two daughter cells with maturity $g(m)$. We assume
that $g:[0,1] \rightarrow [0,1]$ is a continuous and strictly
increasing function, continuously differentiable on $[0,1)$ and
such that $g(m)\leq m$ for $m\in [0,1]$. We also assume, for
technical reason and without loss of generality, that
\begin{displaymath}
\lim_{m\to 1}g^{\prime }(m)=+\infty.
\end{displaymath}
Then we can set
\begin{displaymath}
g^{-1}(m)=1, \qquad \textrm{ for } m>g(1).
\end{displaymath}
This means that the function $g^{-1}:[0,1] \rightarrow [0,1]$ is
continuously differentiable and satisfies $(g^{-1})^{\prime
}(m)=0$, for $m>g(1)$. We set
\begin{displaymath}
\Delta(m)=\Theta(g^{-1}(m)), \qquad \textrm{ for } m\in[0,1].
\end{displaymath}
The quantity $\Delta(m)$ is the maturity of a mother cell at the
point of commitment, when the daughter cells have a maturity $m$
at birth. The function $\Delta : [0,1] \rightarrow [0,1] $ is
continuous and  continuously differentiable on $(0,1]$. It
satisfies $\Delta(0)=0$, $\Delta$ is strictly increasing on
$(0,g(1))$, with $\Theta(m)\leq \Delta(m)$, and
$\Delta(m)=\Theta(1)$ for $m\in[g(1),1]$.

At time $t=0$, the resting and proliferating populations are given
by
\begin{equation}
N(0,m)=\overline{\mu}(m),\label{initialconditionN}
\end{equation}
and
\begin{equation}
P(0,m)=\overline{\Gamma}(m):= \int_0^{\tau(\Theta (m))}
\Gamma(m,a)da,\label{initialconditionP}
\end{equation}
where $\Gamma(m,a)$ is the density of cells with maturity $m$, at
time $t=0$, which have spent a time $a$ in the proliferating
phase, or, equivalently, with age $a$. The functions
$\overline{\mu}$ and $\Gamma$ are supposed to be continuous on
their domains.

We define the sets
\begin{displaymath}
\Omega:=[0,1]\times[0,\tau_{max}],
\end{displaymath}
where $\tau_{max}:=\max_{m\in[0,1]}\tau(m) >0$,
\begin{displaymath}
\Omega_{\Delta} := \Big\{ (m,t)\in\Omega \ ; \ 0\leq t \leq
\tau(\Delta (m)) \Big\},
\end{displaymath}
and
\begin{displaymath}
\Omega_{\Theta} := \Big\{ (m,t)\in\Omega \ ; \ 0\leq t \leq
\tau(\Theta (m)) \Big\}.
\end{displaymath}

\noindent Then, the population densities $N(t,m)$ and $P(t,m)$
satisfy, for $m\in[0,1]$ and $t\geq 0$, the following equations,
\begin{equation}\label{equationN}
\begin{array}{l}
\displaystyle\frac{\partial }{\partial t}N(t,m)+\displaystyle\frac{\partial }{\partial m}(V(m)N(t,m))=-\Big(\delta (m)+\beta \big(m, N(t,m)\big)\Big)N(t,m)\\
\quad \\
+ \left\{ \begin{array}{ll}
\displaystyle 2\xi(t,m)\Gamma\Big( \chi\big(-t,g^{-1}(m)\big),\tau(\Delta(m))-t\Big), &\textrm{if } (m,t)\in\Omega_{\Delta},\\
\quad \\
\displaystyle 2\xi\big(\tau(\Delta(m)),m\big)\beta\Big(\Delta(m),N\big(t-\tau(\Delta(m)),\Delta(m)\big)\Big)\times & \\
\quad \\
\qquad\qquad\qquad\qquad\qquad\qquad
N\big(t-\tau(\Delta(m)),\Delta(m)\big), &\textrm{if }
(m,t)\notin\Omega_{\Delta},
\end{array} \right.
\end{array}
\end{equation}
and
\begin{equation}\label{equationP}
\begin{array}{l}
\displaystyle\frac{\partial }{\partial t}P(t,m)+\displaystyle\frac{\partial }{\partial m}(V(m)P(t,m))=-\gamma(m)P(t,m)+\beta \big(m, N(t,m)\big)N(t,m)\\
\quad \\
-\left\{ \setlength\arraycolsep{2pt} \begin{array}{ll}
\displaystyle \pi(t,m)\Gamma\Big( \chi\big(-t,m\big),\tau(\Theta(m))-t\Big), & \quad \textrm{if } (m,t)\in\Omega_{\Theta},\\
\quad \\
\displaystyle \pi\big(\tau(\Theta(m)),m\big)\beta\Big(\Theta(m),N\big(t-\tau(\Theta(m)),\Theta(m)\big)\Big)\times\\
\quad \\
\quad\qquad\qquad\qquad\qquad\qquad\quad
N\big(t-\tau(\Theta(m)),\Theta(m)\big), & \quad \textrm{if }
(m,t)\notin\Omega_{\Theta},
\end{array} \right.
\end{array}
\end{equation}
where the mappings $\xi:\Omega_{\Delta}\to[0,+\infty)$ and
$\pi:\Omega_{\Theta}\to[0,+\infty)$ are continuous and satisfy
\begin{displaymath}
\xi(\cdot,m)=0 \quad \textrm{ if } \quad m>g(1),
\end{displaymath}
because, from the definition of $g$, a daughter cell can not have
a maturity greater than $g(1)$.

In Equation (\ref{equationN}), the first term in the right hand
side accounts for cellular loss, through cells death ($\delta$)
and introduction in the proliferating phase ($\beta$). The second
term describes the contribution of proliferating cells, one
generation time ago. In a first time, cells can only proceed from
cells initially in the proliferating phase ($\Gamma$). Then, after
one generation time, all cells have divided and the contribution
can only comes from resting cells which have been introduced in
the proliferating phase one generation time ago.

\noindent The factor $2$ always accounts for mitosis. The quantity
$\xi(t,m)$ is for the rate of surviving cells.

\noindent In Equation (\ref{equationP}), the first term in the
right hand side also accounts for cellular loss, whereas the
second term is for the contribution of the resting phase. The
third term describes the same situation as in Equation
(\ref{equationN}), however, in this case, cells leave the
proliferating phase to the resting one. The quantity $\pi(t,m)$ is
also for the rate of surviving cells.

\noindent We can observe two different behaviours of the rates of
surviving cells, in the two phases. In a first time, they depend
on time and maturity, and after a certain time, they only depend
on the maturity variable. When the process of production of blood
cells has just begun, the only cells which divide come from the
initial proliferating phase population. But after one cellular
cycle, that means when $t>\tau(\Delta(m))$ (respectively,
$t>\tau(\Theta(m))$), the amount of cells only comes from resting
cells (respectively, proliferating cells) which have been
introduced in the proliferating phase (respectively, resting
phase) one generation time ago. Consequently, we take into account
the duration of the cell cycle, and not the present time.

\noindent Equations (\ref{equationN}) and (\ref{equationP}) are
derived, after integration, from an age-maturity structured model,
presented by the authors in \cite{adimycrauste}. In fact, the
rates $\xi$ and $\pi$ are explicitly given (see
\cite{adimycrauste}) by
\begin{displaymath}
\xi(t,m)=(g^{-1})^{\prime }(m)\exp\bigg\{ -\int_{0}^{t}\Big(
\gamma\big(\chi(-s,g^{-1}(m))\big)+V^{\prime
}\big(\chi(-s,g^{-1}(m))\big)\Big) ds\bigg\},
\end{displaymath}
and
\begin{displaymath}
\pi(t,m)=\alpha(m)\exp\bigg\{ -\int_{0}^{t}\Big(
\gamma\big(\chi(-s,m)\big)+V^{\prime }\big(\chi(-s,m)\big)\Big)
ds\bigg\},
\end{displaymath}
with $\alpha:[0,1]\to[0,+\infty)$ a positive and continuous
function, such that $\alpha(0)=1$.

\noindent In the following, to simplify the notations, we will
denote by $\overline{\xi}$ and $\overline{\pi}$ the quantities
\begin{displaymath}
\overline{\xi}(m)=\xi\big(\tau(\Delta(m)),m\big),
\end{displaymath}
and
\begin{displaymath}
\overline{\pi}(m)=\pi\big(\tau(\Theta(m)),m\big).
\end{displaymath}

\noindent One can remark that the solutions of Equation
(\ref{equationN}) do not depend on the solutions of Equation
(\ref{equationP}), whereas the converse is not true.

Before we study the asymptotic behaviour of the solutions of
Problem (\ref{initialconditionN})-(\ref{equationP}), we establish
an integrated formulation of this problem. We first extend $N$ by
setting
\begin{equation} \label{prolongementN}
N(t,m)=\overline{\mu}(m), \quad \textrm{ for } t\in[-\tau_{max},
0] \textrm{ and } m\in[0,1].
\end{equation}
One can remark that this extension does not influence the system.

\noindent We also define two mappings,
$F:[0,+\infty)\times[0,1]\times\mathbb{R}\to \mathbb{R}$ and
$G:[0,+\infty)\times[0,1]\times\mathbb{R}\to \mathbb{R}$, by
\begin{equation} \label{F}
F(t,m,x)=\left\{ \begin{array}{ll}
\displaystyle 2\xi(t,m)\Gamma\Big( \chi\big(-t,g^{-1}(m)\big),\tau(\Delta(m))-t\Big), & \textrm{if } (m,t)\in\Omega_{\Delta},\\
\quad \\
\displaystyle 2\overline{\xi}(m)\beta\big(\Delta(m),x\big)x, &
\textrm{if } (m,t)\notin\Omega_{\Delta},
\end{array} \right.
\end{equation}
and
\begin{equation} \label{G}
G(t,m,x)= \left\{ \setlength\arraycolsep{2pt} \begin{array}{ll}
\displaystyle \pi(t,m)\Gamma\Big( \chi\big(-t,m\big),\tau(\Theta(m))-t\Big), & \textrm{if } (m,t)\in\Omega_{\Theta},\\
\quad \\
\displaystyle \overline{\pi}(m)\beta\big(\Theta(m),x\big)x, &
\textrm{if } (m,t)\notin\Omega_{\Theta}.
\end{array} \right.
\end{equation}

\noindent We denote by $C[0,1]$ the space of continuous functions
on $[0,1]$, endowed with the supremum norm $||.||$, defined by
\begin{displaymath}
||v||=\sup_{m\in[0,1]} |v(m)|, \qquad \textrm{ for } v\in C[0,1],
\end{displaymath}
and we consider the unbounded closed linear operator $A :
D(A)\subset C[0,1]\to C[0,1]$ defined by
\begin{displaymath}
D(A)=\Big\{ u\in C[ 0,1] \ ; u \textrm{ differentiable on }(0,1],
u^{\prime }\in C(0,1], \ \lim_{x\rightarrow 0}V(x)u^{\prime }(x)=0
\Big\}
\end{displaymath}
and
\begin{displaymath}
Au(x)=\left\{ \begin{array}{ll}
-(\delta(x)+V^{\prime }(x))u(x)-V(x)u^{\prime }(x), & \quad \textrm{if } x\in(0,1], \\
-(\delta(0)+V^{\prime }(0))u(0), & \quad \textrm{if }x=0.
\end{array}\right.
\end{displaymath}

\begin{proposition} \label{propsemigroupe}
The operator $A$ is the infinitesimal generator of the strongly
continuous semigroup $(T(t))_{t\geq 0}$ defined on $C[0,1]$ by
\begin{displaymath}
(T(t)\psi)(x)=K(t,x)\psi(\chi(-t,x)), \quad \textrm{ for } \psi\in
C[0,1], \ t\geq 0 \textrm{ and } x\in[0,1],
\end{displaymath}
where
\begin{displaymath}
K(t,x)=\exp \left\{ -\int_{0}^{t}\Big(
\delta\big(\chi(-s,x)\big)+V^{\prime }\big(\chi(-s,x)\big)\Big)ds
\right\}.
\end{displaymath}
\end{proposition}

\begin{proof}
The proof is similar to the proof of Proposition 2.4 in Dyson
\emph{et al.} \cite{webb1996}.
\end{proof}

\noindent Now, by using the variation of constants formula
associated to the $C_0$-semigroup $(T(t))_{t\geq 0}$, we can write
an integrated formulation of Problem
(\ref{initialconditionN})-(\ref{equationP}).

\noindent Let $C(\Omega)$ be the space of continuous functions on
$\Omega$, endowed with the norm
\begin{displaymath}
\| \Upsilon \|_{\Omega}:=\sup_{(m,a)\in\Omega}|\Upsilon(m,a)|,
\quad \textrm{ for } \Upsilon\in C(\Omega).
\end{displaymath}

\noindent Let $\overline{\mu}\in C[0,1]$ and $\Gamma\in
C(\Omega)$. An \emph{integrated solution} of Problem
(\ref{initialconditionN})-(\ref{equationP}) is a continuous
solution of the system
\begin{equation} \label{integratedformulationN}
\begin{array}{l}
N(t,m)=K(t,m)\overline{\mu}\big(\chi(-t,m)\big)\\
 \\
-\displaystyle\int_0^t \!\! K(t-s,m)\beta\Big(\chi(-(t-s),m),N\big(s,\chi(-(t-s),m)\big)\Big)N\big(s,\chi(-(t-s),m)\big) ds\\
 \\
+\displaystyle\int_0^t \!\! K(t-s,m)
F\Big(s,\chi(-(t-s),m),N\big(s-\tau(\Delta(\chi(-(t-s),m))),\Delta(\chi(-(t-s),m))\big)\Big)
ds,
\end{array}
\end{equation}
and
\begin{equation} \label{integratedformulationP}
\begin{array}{l}
P(t,m)=H(t,m)\overline{\Gamma}\big(\chi(-t,m)\big)\\
 \\
+\displaystyle\int_0^t \!\! H(t-s,m)\beta\Big(\chi(-(t-s),m),N\big(s,\chi(-(t-s),m)\big)\Big)N\big(s,\chi(-(t-s),m)\big) ds\\
 \\
-\displaystyle\int_0^t \!\! H(t-s,m)
G\Big(s,\chi(-(t-s),m),N\big(s-\tau(\Theta(\chi(-(t-s),m))),\Theta(\chi(-(t-s),m))\big)\Big)
ds,
\end{array}
\end{equation}
for $t\geq 0$ and $m\in[0,1]$, where $F$ and $G$ are given by
(\ref{F}) and (\ref{G}), $\overline{\Gamma}$ is given by
(\ref{initialconditionP}) and
\begin{displaymath}
H(t,m):=\exp \left\{ -\int_{0}^{t}\Big(
\gamma\big(\chi(-s,m)\big)+V^{\prime }\big(\chi(-s,m)\big)\Big)
ds\right\}, \quad \textrm{ for } t\geq 0 \textrm{ and } m\in[0,1].
\label{defksi}
\end{displaymath}

\noindent We can easily prove (see \cite{adimycrauste}), under the
assumptions that the function $x\mapsto \beta(m,x)$ is uniformly
bounded and the function $x\mapsto x\beta(m,x)$ is locally
Lipschitz continuous for all $m\in[0,1]$, that Problem
(\ref{integratedformulationN})-(\ref{integratedformulationP}) has
a unique continuous global solution
$(N^{\overline{\mu},\Gamma},P^{\overline{\mu},\Gamma})$, for
initial conditions $(\overline{\mu},\Gamma)\in C[0,1]\times
C(\Omega)$.

\section{A uniqueness result}\label{sectionuniquenessresult}

In this section, we establish more than uniqueness. Indeed, we
show a result which stresses, for a finite time, the dependence of
the entire population with the small maturity cells population. It
has been shown for the first time by Dyson \emph{et al.}
\cite{webb1996}, for a model with a constant delay. We will see
that this result is important in order to obtain the asymptotic
behaviour of the solutions of
(\ref{integratedformulationN})-(\ref{integratedformulationP}).

\noindent We first assume that
\begin{equation}
\Delta (m)<m, \qquad \textrm{ for all }
m\in(0,1].\label{inegtauvarianlin}
\end{equation}
This condition is equivalent to
\begin{equation}
\tau(\Delta(m))>\int_{m}^{g^{-1}(m)}\frac{ds}{V(s)},\quad \textrm{
for } m\in(0,1].  \label{inegtauvarianlin2}
\end{equation}
This equivalence is immediate when one notices that, from
(\ref{propertyoftheta}),
\begin{displaymath}
\Delta(m)=\chi\big(-\tau(\Delta
(m)),g^{-1}(m)\big)=h^{-1}\big(h(g^{-1}(m))e^{-\tau(\Delta
(m))}\big).
\end{displaymath}
Since the quantity $\int_{m}^{g^{-1}(m)}\frac{ds}{V(s)}$
represents the time required for a cell with maturity $m$, at
birth, to reach the maturity of its mother at the cytokinesis
point (the point of division), Condition (\ref{inegtauvarianlin2})
means that, in the proliferating phase, cells have enough time to
reach the maturity of their mother.

\noindent Condition (\ref{inegtauvarianlin}) implies in particular
that
\begin{displaymath}
\Theta(1):=\Delta(g(1))<g(1).
\end{displaymath}

\noindent From now on, and throughout this section, we assume that
the function $x\mapsto \beta(m,x)$ is uniformly bounded, the
function $x\mapsto x\beta(m,x)$ is locally Lipschitz continuous
for all $m\in[0,1]$, and that Condition (\ref{inegtauvarianlin})
holds.

\noindent For $b\in(0,1]$ and $\psi\in C[0,1]$, we define $\| .
\|_b$ as follows
\begin{displaymath}
\|\psi\|_b := \sup_{m\in[0,b]}|\psi(m)|.
\end{displaymath}

\noindent We first show the following proposition.

\begin{proposition}\label{propositionuniqueness}
Let $\overline{\mu}_1, \overline{\mu}_2\in C[0,1]$ and $\Gamma_1,
\Gamma_2 \in C(\Omega)$. If there exists $0<b<1$ such that
\begin{equation}\label{hypN}
\overline{\mu}_{1}(m)=\overline{\mu}_{2}(m) \quad \textrm{and}
\quad \Gamma_1(m,a)=\Gamma_2(m,a),
\end{equation}
for $m\in[0,b]$ and $a\in[0,\tau_{max}]$, then,
\begin{equation}\label{propertyN}
N^{\overline{\mu}_1,\Gamma_1}(t,m)=N^{\overline{\mu}_2,\Gamma_2}(t,m),
\quad \textrm{ for } t\geq0 \textrm{ and } m\in[0,g(b)].
\end{equation}
\end{proposition}

\begin{proof}
We suppose that there exists $b\in(0,1)$ such that (\ref{hypN})
holds. Let $T>0$ be given, and let $t\in(0,T]$ and $m\in[0, g(b)]$
be fixed. Since $h$ is increasing, it follows from (\ref{charach})
that
\begin{displaymath}
\chi(-t,m)\leq m \leq g(b) \leq b.
\end{displaymath}
Then
\begin{displaymath}
\overline{\mu}_1(\chi(-t,m))=\overline{\mu}_2(\chi(-t,m)).
\end{displaymath}
Let $s\in[0,t]$. Since $g^{-1}$ is increasing, then
\begin{displaymath}
\chi\big(-s,g^{-1}\big(\chi(-(t-s),m)\big)\big) \leq
g^{-1}\big(\chi(-(t-s),m)\big) \leq g^{-1}(m) \leq b.
\end{displaymath}
Moreover, if $0\leq s\leq \tau\big(\Delta(\chi(-(t-s),m))\big)$,
then
\begin{displaymath}
\tau\big(\Delta(\chi(-(t-s),m))\big)-s \in [0,\tau_{max}].
\end{displaymath}
Thus, we have
\begin{displaymath}
\begin{array}{l}
\Gamma_1\Big(\chi\big(-s,g^{-1}\big(\chi(-(t-s),m)\big)\big),\tau\big(\Delta(\chi(-(t-s),m))\big)-s\Big)\\
\quad\\
=
\Gamma_2\Big(\chi\big(-s,g^{-1}\big(\chi(-(t-s),m)\big)\big),\tau\big(\Delta(\chi(-(t-s),m))\big)-s\Big).
\end{array}
\end{displaymath}
Since the solutions $N^{\overline{\mu}_1,\Gamma_1}(t,m)$ and
$N^{\overline{\mu}_2,\Gamma_2}(t,m)$ of Equation
(\ref{integratedformulationN}) are continuous and satisfy
\begin{displaymath}
N^{\overline{\mu}_1,\Gamma_1}(0,m)=N^{\overline{\mu}_2,\Gamma_2}(0,m),
\qquad \textrm{ for } m\in[0,b],
\end{displaymath}
then, by using the locally Lipshitz continuous property of the
function $x\mapsto x\beta(m,x)$, we can write
\begin{displaymath}
\begin{array}{l}
|N^{\overline{\mu}_1,\Gamma_1}(t,m)-N^{\overline{\mu}_2,\Gamma_2}(t,m)|\\
\quad\\
\leq \widetilde{K}L \displaystyle\int_0^t |N^{\overline{\mu}_1,\Gamma_1}(s,\chi(-(t-s),m))-N^{\overline{\mu}_2,\Gamma_2}(s,\chi(-(t-s),m))| ds\\
\quad\\
+2\widetilde{K}L\|\overline{\xi}\| \displaystyle\int_0^t |N^{\overline{\mu}_1,\Gamma_1}(s-\tau(\Delta(\chi(-(t-s),m))),\Delta(\chi(-(t-s),m)))\\
\quad\\
\qquad\qquad\qquad -N^{\overline{\mu}_2,\Gamma_2}(s-\tau(\Delta(\chi(-(t-s),m))),\Delta(\chi(-(t-s),m)))| ds,\\
\quad\\
\leq \widetilde{K}L \displaystyle\int_0^t \| N^{\overline{\mu}_1,\Gamma_1}(s,.)-N^{\overline{\mu}_2,\Gamma_2}(s,.) \|_{g(b)} ds\\
\quad\\
+2\widetilde{K}L\|\overline{\xi}\| \displaystyle\int_0^t \|
N^{\overline{\mu}_1,\Gamma_1}(s-\tau(\Delta(\chi(-(t-s),m))),.)-N^{\overline{\mu}_2,\Gamma_2}(s-\tau(\Delta(\chi(-(t-s),m))),.)
\|_{g(b)} ds,
\end{array}
\end{displaymath}
for $T>0$ small enough, where $L$ is a Lipschitz constant of the
function $x\mapsto x\beta(m,x)$ and $\widetilde{K}$ is defined by
\begin{displaymath}
K(s,m)\leq \widetilde{K}, \quad \textrm{ for } s\in[0,T] \textrm{
and } m\in[0,1].
\end{displaymath}
The extension given by (\ref{prolongementN}) allows to give sense
to the integral terms in the above inequality.

\noindent Let $\theta\in[-\tau_{max},0]$ be given. If $t+\theta <
0$, then
$N^{\overline{\mu}_1,\Gamma_1}(t+\theta,m)=N^{\overline{\mu}_2,\Gamma_2}(t+\theta,m)$.
If $t+\theta \geq 0$, then
\begin{displaymath}
\begin{array}{l}
|N^{\overline{\mu}_1,\Gamma_1}(t+\theta,m)-N^{\overline{\mu}_2,\Gamma_2}(t+\theta,m)|\\
\quad\\
\leq \widetilde{K}L \displaystyle\int_0^{t+\theta} \| N^{\overline{\mu}_1,\Gamma_1}(s,.)-N^{\overline{\mu}_2,\Gamma_2}(s,.) \|_{g(b)} ds\\
\quad\\
+2\widetilde{K}L\|\overline{\xi}\| \displaystyle\int_0^{t+\theta} \| N^{\overline{\mu}_1,\Gamma_1}(s-\tau(\Delta(\chi(-(t+\theta-s),m))),.)-N^{\overline{\mu}_2,\Gamma_2}(s-\tau(\Delta(\chi(-(t+\theta-s),m))),.) \|_{g(b)} ds,\\
\quad\\
\leq \widetilde{K}L(1+2\|\overline{\xi}\|)\displaystyle\int_0^t
\sup_{\overline{\theta}\in[-\tau_{max},0]}\|
N^{\overline{\mu}_1,\Gamma_1}(s+\overline{\theta},.)-N^{\overline{\mu}_2,\Gamma_2}(s+\overline{\theta},.)
\|_{g(b)} ds.
\end{array}
\end{displaymath}
It follows that
\begin{displaymath}
\setlength\arraycolsep{2pt}
\begin{array}{l}
\displaystyle\sup_{\theta\in[-\tau_{max},0]}\| N^{\overline{\mu}_1,\Gamma_1}(t+\theta,.)-N^{\overline{\mu}_2,\Gamma_2}(t+\theta,.) \|_{g(b)} \\
\qquad\qquad\qquad\qquad\leq
\widetilde{K}L(1+2\|\overline{\xi}\|)\displaystyle\int_0^t
\displaystyle\sup_{\theta\in[-\tau_{max},0]}\|
N^{\overline{\mu}_1,\Gamma_1}(s+\theta,.)-N^{\overline{\mu}_2,\Gamma_2}(s+\theta,.)
\|_{g(b)} ds.
\end{array}
\end{displaymath}
By using the Gronwall's Inequality, we obtain
\begin{displaymath}
\sup_{\theta\in[-\tau_{max},0]}\|
N^{\overline{\mu}_1,\Gamma_1}(t+\theta,.)-N^{\overline{\mu}_2,\Gamma_2}(t+\theta,.)
\|_{g(b)}=0.
\end{displaymath}
In particular,
\begin{displaymath}
\|
N^{\overline{\mu}_1,\Gamma_1}(t,.)-N^{\overline{\mu}_2,\Gamma_2}(t,.)
\|_{g(b)}=0, \qquad \textrm{ for } t\in(0,T].
\end{displaymath}
By steps, this result holds for all $T>0$, therefore
(\ref{propertyN}) is satisfied and the proof is complete.
\end{proof}

\noindent Now, let $0<\overline{b}<g(1)$ be fixed and consider the
sequence $(b_n)_{n\in\mathbb{N}}$ defined by
\begin{equation} \label{suitebn}
b_{0}=\overline{b} \quad \textrm{ and } \quad b_{n+1}=\left\{
\begin{array}{ll}
\Delta^{-1}(b_{n}), & \textrm{ if } b_{n}\in[0,\Theta(1)),\\
            &   \\
g(1),           & \textrm{ if } b_{n}\in[\Theta(1), g(1)].
\end{array}\right.
\end{equation}
The sequence $(b_{n})_{n\in\mathbb{N}}$ represents the
transmission of the maturity between two successive generations,
$n$ and $n+1$. The following result is immediate.

\begin{lemma}\label{lemma2}
If $0<\overline{b}<\Theta(1):=\Delta(g(1))$, then there exists
$N\in\mathbb{N}$ such that $b_{N}<\Theta(1)\leq b_{N+1}\leq g(1)$.
\end{lemma}

\noindent We give now a first result, which emphasizes the strong
link between the process of production of cells and the population
of stem cells. A similar result has been proved by Adimy and
Pujo-Menjouet \cite{adimypujo2003} in the linear case.

\begin{theorem}\label{smallmaturitycells}
Let $\overline{\mu}_1, \overline{\mu}_2\in C[0,1]$ and $\Gamma_1,
\Gamma_2 \in C(\Omega)$. If there exists $0<b<1$ such that
(\ref{hypN}) holds, then, there exists $\overline{t}>0$ such that
\begin{displaymath}
N^{\overline{\mu}_1,\Gamma_1}(t,m)=N^{\overline{\mu}_2,\Gamma_2}(t,m),
\end{displaymath}
for $m\in[0,g(1)]$ and $t\geq\overline{t}$, where $\overline{t}$
can be chosen to be
\begin{equation}\label{tbar}
\overline{t}=\ln \Bigg[\frac{h(g(1))}{h(g(b))}\Bigg]
+(N+2)\tau_{max},
\end{equation}
and $N\in \mathbb{N}$ is given by Lemma \ref{lemma2}, for
$\overline{b}=g(b)$. Furthermore,
\begin{displaymath}
N^{\overline{\mu}_1,\Gamma_1}(t,m)=N^{\overline{\mu}_2,\Gamma_2}(t,m),
\end{displaymath}
for $m\in [g(1),1]$ and $t\geq
\overline{t}+\tau_{max}-\ln\big(h(g(1))\big)=(N+3)\tau_{max}-\ln\big(h(g(b))\big)$.
\end{theorem}

\begin{proof}
Let $\overline{b}=g(b)$. Since $g$ is increasing, then
$\overline{b}<g(1)$. Proposition \ref{propositionuniqueness}
implies that
\begin{displaymath}
N^{\overline{\mu}_1,\Gamma_1}(t,m)=N^{\overline{\mu}_2,\Gamma_2}(t,m),
\quad \textrm{ for } t\geq0 \textrm{ and } m\in[0,\overline{b}].
\end{displaymath}
Let us reconsider the sequence $(b_{n})_{n\in\mathbb{N}}$, given
by (\ref{suitebn}), and let us consider the sequence
$(t_{n})_{n\in \mathbb{N}}$ defined by
\begin{equation} \label{recurr1}
\left\{ \begin{array}{l}
t_{n+1}=t_{n}+\ln \Bigg[\displaystyle\frac{h(b_{n+1})}{h(b_{n})}\Bigg]+\tau_{max},\\
t_{0}=0.
\end{array} \right.
\end{equation}
Then,
\begin{displaymath}
t_{n}=\ln \Bigg[ \frac{h(b_{n})}{h(g(b))}\Bigg] +n\tau_{max}.
\end{displaymath}
The sequence $(b_{n})_{n\in \mathbb{N}}$ is increasing. Then, the
sequence $(t_{n})_{n\in \mathbb{N}}$ is also increasing.  We are
going to prove, by induction, the following result
\begin{displaymath}
\left( H_{n}\right) : \qquad
N^{\overline{\mu}_1,\Gamma_1}(t,m)=N^{\overline{\mu}_2,\Gamma_2}(t,m),
\quad \textrm{ for } t\geq t_{n} \textrm{ and } m\in[0,b_{n}].
\end{displaymath}
First, $(H_{0})$ is true, from Proposition
\ref{propositionuniqueness}.

\noindent Let suppose that $(H_{n}) $ is true for $n\in
\mathbb{N}$. Let $t\geq t_{n+1}$ and $m\in [0,b_{n+1}]$. Then,
from (\ref{recurr1}),
\begin{displaymath}
t_{n+1}\geq t_{n}+\tau_{max}\geq\tau_{max}.
\end{displaymath}
Since Equation (\ref{integratedformulationN}) is autonomous, its
solutions can be reformulated, for $t\geq t_{n+1}$, as follows
\begin{displaymath}
\setlength\arraycolsep{2pt}
\begin{array}{l}
N^{\overline{\mu}_i,\Gamma_i}(t,m)= K(t-t_{n}-\tau_{max},m)N^{\overline{\mu}_i,\Gamma_i}\Big(t_{n}+\tau_{max},\chi\big(-(t-t_{n}-\tau_{max}),m\big)\Big)\\
 \\
-\displaystyle\int_{t_{n}+\tau_{max}}^{t} K(t-s,m)\beta\bigg(\chi\big(-(t-s),m\big),N^{\overline{\mu}_i,\Gamma_i}\Big(s,\chi\big(-(t-s),m\big)\Big)\bigg)N^{\overline{\mu}_i,\Gamma_i}\Big(s,\chi\big(-(t-s),m\big)\Big)ds\\
 \\
+2\displaystyle\int_{t_{n}+\tau_{max}}^{t} K(t-s,m)\overline{\xi}\Big(\chi\big(-(t-s),m\big)\Big)N^{\overline{\mu}_i,\Gamma_i}\Big(s-\tau\big(\Delta(\chi(-(t-s),m))\big),\Delta(\chi(-(t-s),m))\Big)\times\\
 \\
\qquad\qquad\qquad\beta\bigg(\Delta(\chi(-(t-s),m)),N^{\overline{\mu}_i,\Gamma_i}\Big(s-\tau\big(\Delta(\chi(-(t-s),m))\big),\Delta(\chi(-(t-s),m))\Big)\bigg)
\ ds,
\end{array}
\end{displaymath}
for $i=1,2$. Remark that, from (\ref{charach}),
\begin{displaymath}
\chi\big(-(t-t_{n}-\tau_{max}),m\big)=h^{-1}\Big(h(m)e^{-(t-t_{n}-\tau_{max})}\Big),
\end{displaymath}
and from (\ref{recurr1}),
\begin{displaymath}
\begin{array}{ll}
e^{-(t-t_{n}-\tau_{max})} & =\displaystyle\frac{h(b_{n})}{h(b_{n+1})}e^{-(t-t_{n+1})}, \\
              & \leq \displaystyle\frac{h(b_{n})}{h(b_{n+1})}.
\end{array}
\end{displaymath}
Then, we deduce that
\begin{displaymath}
\begin{array}{rl}
\chi\big(-(t-t_{n}-\tau_{max}),m\big) & \leq h^{-1}\bigg(h(m)\displaystyle\frac{h(b_{n})}{h(b_{n+1})}\bigg), \\
                & \leq h^{-1}\bigg(h(b_{n+1})\displaystyle\frac{h(b_{n})}{h(b_{n+1})}\bigg)=b_{n}.
\end{array}
\end{displaymath}
Hence, $(H_{n})$ implies
\begin{displaymath}
N^{\overline{\mu}_1,\Gamma_1}\Big(t_{n}+\tau_{max},\chi\big(-(t-t_{n}-\tau_{max}),m\big)\Big)=N^{\overline{\mu}_2,\Gamma_2}\Big(t_{n}+\tau_{max},\chi\big(-(t-t_{n}-\tau_{max}),m\big)\Big).
\end{displaymath}
Furthermore, for $t_{n}+\tau_{max}\leq s\leq t$, we have
\begin{displaymath}
s-\tau\big(\Delta(\chi(-(t-s),m))\big)\geq
(t_{n}+\tau_{max})-\tau\big(\Delta(\chi(-(t-s),m)\big)\geq t_{n},
\end{displaymath}
and
\begin{displaymath}
\Delta \big(\chi(-(t-s),m)\big)\leq \Delta (m)\leq \Delta
(b_{n+1})=b_{n}.
\end{displaymath}
Consequently,
\begin{displaymath}
N^{\overline{\mu}_1,\Gamma_1}\Big(s-\tau\big(\Delta(\chi(-(t-s),m)\big),\Delta
(\chi(-(t-s),m))\Big)=N^{\overline{\mu}_2,\Gamma_2}\Big(s-\tau\big(\Delta(\chi(-(t-s),m)),\Delta
(\chi(-(t-s),m))\big)\Big).
\end{displaymath}
Then, we obtain that \begin{displaymath}
\setlength\arraycolsep{2pt}
\begin{array}{l}
|N^{\overline{\mu}_1,\Gamma_1}(t,m)-N^{\overline{\mu}_2,\Gamma_2}(t,m)|\leq\\
 \\
\displaystyle\int_{t_{n}+\tau_{max}}^{t}\!\!\!\!K(t-s,m)\bigg| \beta\Big(\chi(-(t-s),m),N^{\overline{\mu}_1,\Gamma_1}\big(s,\chi(-(t-s),m)\big)\Big)N^{\overline{\mu}_1,\Gamma_1}\big(s,\chi(-(t-s),m)\big)\\
\\
\quad\quad\quad\quad\quad\quad\quad
-\beta\Big(\chi(-(t-s),m),N^{\overline{\mu}_2,\Gamma_2}\big(s,\chi(-(t-s),m)\big)\Big)N^{\overline{\mu}_2,\Gamma_2}\big(s,\chi(-(t-s),m)\big)\bigg|ds,
\end{array}
\end{displaymath}
and, by using the Gronwall's inequality, we deduce that
$(H_{n+1})$ is true. Consequently, $(H_{n})$ is true for
$n\in\mathbb{N}$.

\noindent In particular, $(H_{n})$ holds for $n=N+2$, where $N$ is
given by Lemma \ref{lemma2}, with $\overline{b}=g(b)$. Since
$b_{N+1}\in[\Theta(1),g(1)]$, then $b_{N+2}=g(1)$. We deduce that
\begin{equation}\label{hypN2}
N^{\overline{\mu}_1,\Gamma_1}(t,m)=N^{\overline{\mu}_2,\Gamma_2}(t,m)
\quad \textrm{ for }m\in[0,g(1)] \textrm{ and }t\geq\overline{t},
\end{equation}
where $\overline{t}=t_{N+2}$ is given by (\ref{tbar}).

\noindent Finally, take $m\in[g(1),1]$ and
$t\geq\overline{t}+\tau_{max}$. We can write, for $i=1,2$,
\begin{displaymath}
\setlength\arraycolsep{2pt}
\begin{array}{l}
N^{\overline{\mu}_i,\Gamma_i}(t,m)=K(t-\overline{t}-\tau_{max},m)N^{\overline{\mu}_i,\Gamma_i}(\overline{t}+\tau_{max},\chi(-(t-\overline{t}-\tau_{max}),m))\\
 \\
-\displaystyle\int_{\overline{t}+\tau_{max}}^{t} K(t-s,m)\beta\bigg(\chi\big(-(t-s),m\big),N^{\overline{\mu}_i,\Gamma_i}\Big(s,\chi\big(-(t-s),m\big)\Big)\bigg)N^{\overline{\mu}_i,\Gamma_i}\Big(s,\chi\big(-(t-s),m\big)\Big)ds\\
 \\
+2\displaystyle\int_{\overline{t}+\tau_{max}}^{t} K(t-s,m)\overline{\xi}\Big(\chi\big(-(t-s),m\big)\Big)N^{\overline{\mu}_i,\Gamma_i}\Big(s-\tau\big(\Delta(\chi(-(t-s),m))\big),\Delta(\chi(-(t-s),m))\Big)\times\\
 \\
\quad\qquad\qquad\qquad\beta\bigg(\Delta(\chi(-(t-s),m)),N^{\overline{\mu}_i,\Gamma_i}\Big(s-\tau\big(\Delta(\chi(-(t-s),m))\big),\Delta(\chi(-(t-s),m))\Big)\bigg)
\ ds.
\end{array}
\end{displaymath}
Let $\overline{t}+\tau_{max} \leq s\leq t$. Then,
\begin{equation*}
s-\tau\big(\Delta(\chi(-(t-s),m))\big)\geq
(\overline{t}+\tau_{max})-\tau\big(\Delta(\chi(-(t-s),m)\big)\geq
\overline{t}.
\end{equation*}
Consequently, if
\begin{equation*}
\chi(-(t-s),m)\leq g(1),
\end{equation*}
then,
\begin{equation*}
\Delta(\chi(-(t-s),m))\leq \Delta(g(1)) < g(1),
\end{equation*}
and (\ref{hypN2}) implies that
\begin{equation*}
N^{\overline{\mu}_1,\Gamma_1}\Big(s-\tau\big(\Delta
(\chi(-(t-s),m))\big),\Delta(\chi(-(t-s),m))\Big)=N^{\overline{\mu}_2,\Gamma_2}\Big(s-\tau\big(\Delta
(\chi(-(t-s),m))\big),\Delta(\chi(-(t-s),m))\Big).
\end{equation*}
On the other hand, if
\begin{equation*}
\chi\big(-(t-s),m\big)>g(1),
\end{equation*}
then, from the definition of $\overline{\xi}$, we have
\begin{equation*}
\overline{\xi}\Big(\chi\big(-(t-s),m\big)\Big)=0.
\end{equation*}
Furthermore, by remarking that $\ln (h(m))\leq 0$, for all $m\in
(0,1]$, then we deduce, for $m\in[g(1),1]$ and $t\geq
\overline{t}+\tau_{max}-\ln (h(g(1)))$, that
\begin{equation*}
\begin{array}{lll}
\chi\big(-(t-\overline{t}-\tau_{max}),m\big)&=& h^{-1}\Big(h(m)e^{-(t-\overline{t}-\tau_{max})}\Big), \\
                         & &  \\
                         &\leq& h^{-1}\big(h(m)h(g(1))\big),\\
                     &   &  \\
                     &\leq& h^{-1}\big(h(g(1))\big)=g(1). \\
\end{array}
\end{equation*}
Hence,
\begin{equation*}
N^{\overline{\mu}_1,\Gamma_1}\Big(\overline{t}+\tau_{max},\chi\big(-(t-\overline{t}-\tau_{max}),m\big)\Big)=N^{\overline{\mu}_2,\Gamma_2}\Big(\overline{t}+\tau_{max},\chi\big(-(t-\overline{t}-\tau_{max}),m\big)\Big).
\end{equation*}
Using once again the Gronwall's inequality, we conclude that
 \begin{equation*}
N^{\overline{\mu}_1,\Gamma_1}(t,m)=N^{\overline{\mu}_2,\Gamma_2}(t,m),
\quad \textrm{ for } m\in [g(1),1] \textrm{ and } t\geq
\overline{t}+\tau_{max}-\ln (h(g(1))).
\end{equation*}
This completes the proof.
\end{proof}

\begin{corollary} \label{coraplascticanemia}
Let $\overline{\mu}\in C[0,1]$ and $\Gamma\in C(\Omega)$. If there
exists $0<b<1$ such that
\begin{displaymath}
\overline{\mu}(m)=0 \quad \textrm{and} \quad \Gamma(m,a)=0, \quad
\textrm{ for } m\in[0,b] \textrm{  and } a\in[0,\tau_{max}],
\end{displaymath}
then,
\begin{displaymath}
N^{\overline{\mu},\Gamma}(t,m)=0, \quad \textrm{ for } m\in [0,1]
\textrm{ and } t\geq (N+3)\tau_{max}-\ln\big(h(g(b))\big),
\end{displaymath}
where $N\in \mathbb{N}$ is given by Lemma \ref{lemma2}, for
$\overline{b}=g(b)$.
\end{corollary}

This result stresses the dependence of the production of cells
with the population of stem cells. In particular, if the stem
cells population is defective in the initial stage, then the
entire population is doomed to extinction in a finite time. This
situation describes what usually happens with the aplastic anemia,
a disease which yields to injury or destruction of pluripotential
stem cells.

In the next corollary, we show that the proliferating population
depends also strongly on the stem cells population.

\begin{corollary}\label{smallmaturityprol}
Let $\overline{\mu}_1, \overline{\mu}_2\in C[0,1]$ and $\Gamma_1,
\Gamma_2 \in C(\Omega)$. If there exists $0<b<1$ such that
(\ref{hypN}) holds, then
\begin{displaymath}
P^{\overline{\mu}_1,\Gamma_1}(t,m)=P^{\overline{\mu}_2,\Gamma_2}(t,m),
\qquad \textrm{ for } m\in [0,1] \textrm{ and } t\geq
(N+3)\tau_{max}-\ln\big(h(g(b))\big),
\end{displaymath}
where $N\in \mathbb{N}$ is given by Lemma \ref{lemma2}, for
$\overline{b}=g(b)$.
\end{corollary}

\begin{proof}
The proof is immediate by using Theorem \ref{smallmaturitycells},
Equation (\ref{integratedformulationP}) and a method of steps.
\end{proof}

\section{Behaviour of the immature cells population} \label{sectionbehaviourimmaturecells}

In this section, we investigate the behaviour of the immature
cells population, that means, the population of cells with
maturity $m=0$.

Let $\overline{\mu}\in C[0,1]$ and $\Gamma\in C(\Omega)$ be fixed.
Let us consider the continuous solutions
$N^{\overline{\mu},\Gamma}(t,m)$ and
$P^{\overline{\mu},\Gamma}(t,m)$ of Problem
(\ref{integratedformulationN})-(\ref{integratedformulationP}). We
set $x(t)=N^{\overline{\mu},\Gamma}(t,0)$ and
$y(t)=P^{\overline{\mu},\Gamma}(t,0)$, for all $t\geq0$. Then,
$\big(x(t),y(t)\big)$ is solution of the system
\begin{equation}\label{equationmatzerox1}
\begin{array}{rcl}
x(t)\!\!\!&=&\!\!\!e^{-\rho t}\overline{\mu}(0) - \displaystyle\int_0^t e^{-\rho(t-s)}\beta(0,x(s))x(s) ds\\
    \!\!\!&+&\!\!\!\left\{\begin{array}{ll} 2\displaystyle\int_0^t e^{-\rho(t-s)}\xi(s,0)\Gamma(0,r-s) ds,&\!\!\!\!\!\! \textrm{ for } t\in[0,r],\\ 2\displaystyle\int_0^r e^{-\rho(t-s)}\xi(s,0)\Gamma(0,r-s) ds+2\overline{\xi}(0)\displaystyle\int_r^t e^{-\rho(t-s)}\beta(0,x(s-r))x(s-r) ds,&\!\!\!\!\!\! \textrm{ for } r\leq t,\end{array}\right.
\end{array}
\end{equation}
and
\begin{equation}\label{equationmatzeroy1}
\begin{array}{rcl}
y(t)\!\!\!&=&\!\!\!e^{-\eta t}\overline{\Gamma}(0)+\displaystyle\int_0^t e^{-\eta(t-s)}\beta(0,x(s))x(s) ds\\
    \!\!\!&-&\!\!\!\left\{\begin{array}{ll} \displaystyle\int_0^t e^{-\eta(t-s)}\pi(s,0)\Gamma(0,r-s) ds,&\!\!\!\!\!\!\textrm{ for } t\in[0,r],\\ \displaystyle\int_0^r e^{-\eta(t-s)}\pi(s,0)\Gamma(0,r-s) ds+\overline{\pi}(0)\displaystyle\int_r^t e^{-\eta(t-s)}\beta(0,x(s-r))x(s-r) ds,&\!\!\!\!\!\! \textrm{ for } r\leq t,\end{array}\right.
\end{array}
\end{equation}
where $\rho:=\delta(0)+V^{\prime}(0)$,
$\eta:=\gamma(0)+V^{\prime}(0)$ and $r:=\tau(0)>0$.

\noindent Let us recall that $\overline{\xi}(0)=\xi(r,0)$,
$\overline{\pi}(0)=\pi(r,0)$ and
\begin{displaymath}
\overline{\Gamma}(0)=\int_0^r \Gamma(0,a)da.
\end{displaymath}
Then, we easily deduce that System
(\ref{equationmatzerox1})-(\ref{equationmatzeroy1}) is equivalent
to the system
\begin{equation}\label{equationmatzero2}
\left\{\begin{array}{rcl}
\displaystyle\frac{dx}{dt}(t)&=&-\big(\rho+\beta(0,x(t))\big)x(t)+\left\{\begin{array}{ll} 2\xi(t,0)\Gamma(0,r-t),&\textrm{ for } 0\leq t\leq r,\\  2\overline{\xi}(0)\beta(0,x(t-r))x(t-r),&\textrm{ for }r\leq t,\end{array}\right.\\
x(0)&=&\overline{\mu}(0),
\end{array}\right.
\end{equation}
\begin{equation}\label{equationmatzeroy}
\left\{\begin{array}{rcl}
\displaystyle\frac{dy}{dt}(t)&=&-\eta y(t)+\beta(0,x(t))x(t)-\left\{\begin{array}{ll}\pi(t,0)\Gamma(0,r-t),&\textrm{ for } 0\leq t\leq r, \\ \overline{\pi}(0)\beta(0,x(t-r))x(t-r),&\textrm{ for }r\leq t,\end{array}\right.\\
y(0)&=&\overline{\Gamma}(0).
\end{array}\right.
\end{equation}
Of course, at $t=r$, the derivatives in (\ref{equationmatzero2})
and (\ref{equationmatzeroy}) represent the right-hand side and the
left-hand side derivatives.

\noindent First, consider the system, for $t\in[0,r]$,
\begin{equation}\label{equationmatzerophipsi}
\left\{\begin{array}{rcl}
\vspace{1ex}\displaystyle\frac{d\phi}{dt}(t)&=&-\big(\rho+\beta(0,\phi(t))\big)\phi(t)+2\xi(t,0)\Gamma(0,r-t),\\
\displaystyle\frac{d\psi}{dt}(t)&=&-\eta
\psi(t)+\beta(0,\phi(t))\phi(t)-\pi(t,0)\Gamma(0,r-t),
\end{array}\right.
\end{equation}
with
\begin{equation}\label{equationmatzerophipsiic}
\left\{\begin{array}{rcl}
\phi(0)&=&\overline{\mu}(0),\\
\psi(0)&=&\overline{\Gamma}(0).
\end{array}\right.
\end{equation}
It is obvious that, under the assumptions that the function
$x\mapsto \beta(0,x)$ is bounded and the function $x\mapsto
x\beta(0,x)$ is locally Lipschitz continuous, Problem
(\ref{equationmatzerophipsi})-(\ref{equationmatzerophipsiic}) has
a unique solution $\big(\phi(t),\psi(t)\big)$, for $t\in[0,r]$.
Remark that $\psi(t)$ is explicitly given by
\begin{equation}\label{psiexplicite}
\psi(t)=e^{-\eta t}\int_0^{r-t} \Gamma(0,a) da + \int_0^t
e^{-\eta(t-s)}\beta(0,\phi(s))\phi(s) ds, \qquad \textrm{ for }
t\in[0,r].
\end{equation}
Moreover, if $\overline{\mu}(0)\geq0$ and $\Gamma(0,\cdot)\geq0$,
then $\phi(t)$ and $\psi(t)$ are nonnegative.

\noindent Hence, for $t\geq r$, Problem
(\ref{equationmatzero2})-(\ref{equationmatzeroy}) reduces to the
delay differential system
\begin{equation}\label{equationmatzero}
\frac{dx}{dt}(t)=-\big(\rho+\beta(0,x(t))\big)x(t)+2\overline{\xi}(0)\beta(0,x(t-r))x(t-r),
\end{equation}
\begin{equation}\label{equationmatzeroprime}
\frac{dy}{dt}(t)=-\eta
y(t)+\beta(0,x(t))x(t)-\overline{\pi}(0)\beta(0,x(t-r))x(t-r),
\end{equation}
with, for $t\in[0,r]$,
\begin{equation}\label{equationmatzeroic}
\left\{\begin{array}{rcl}
x(t)&=&\phi(t),\\
y(t)&=&\psi(t).
\end{array}\right.
\end{equation}

\noindent As $\psi(t)$, for $t\in[0,r]$, $y(t)$ is explicitly
given, for $t\geq r$, by
\begin{equation}\label{yexplicite}
y(t)=\int_{t-r}^t e^{-\eta(t-s)}\beta(0,x(s))x(s)ds.
\end{equation}

\begin{proposition}
Assume that the function $x\mapsto \beta(0,x)$ is bounded and the
function $x\mapsto x\beta(0,x)$ is locally Lipshitz continuous.
Let $\overline{\mu}\in C[0,1]$ and $\Gamma\in C(\Omega)$ be given.
Then, Problem (\ref{equationmatzero})-(\ref{equationmatzeroic})
has a unique solution $\big(x^{\phi}(t),y^{\psi}(t)\big)$, defined
for $t\geq0$, where $\big(\phi(t),\psi(t)\big)$ is the unique
solution of
(\ref{equationmatzerophipsi})-(\ref{equationmatzerophipsiic}).
Moreover, $\big(x^{\phi}(t),y^{\psi}(t)\big)$ has a continuous
derivative at $t=r$ if and only if
\begin{equation}\label{condcompatibilite}
\Gamma(0,0)=\beta\big(0,\overline{\mu}(0)\big)\overline{\mu}(0).
\end{equation}
Furthermore, if $\overline{\mu}(0)\geq0$ and
$\Gamma(0,\cdot)\geq0$, then $x^{\phi}(t)$ and $y^{\psi}(t)$ are
nonnegative.
\end{proposition}

\begin{proof}
Existence, uniqueness and regularity of solutions of Problem
(\ref{equationmatzero})-(\ref{equationmatzeroic}) come from Hale
and Verduyn Lunel \cite{halelunel}. The positivity of
$x^{\phi}(t)$ is easily obtained by steps. Moreover, by using
(\ref{yexplicite}), we deduce immediatly the positivity of
$y^{\psi}(t)$.
\end{proof}

In the sequel, we will consider, for a biological reason, only
nonnegative solutions of Problem
(\ref{equationmatzero})-(\ref{equationmatzeroprime}).

\begin{lemma}\label{lemmastability}
If $\lim_{t \to +\infty} x^{\phi}(t)=C$ exists, then $\lim_{t \to
+\infty} y^{\psi}(t)$ exists and is equal to
\begin{displaymath}
\left\{ \begin{array}{ll}\displaystyle\frac{1}{\eta}(1-e^{-\eta
r})\beta(0,C)C,&\qquad \textrm{ if } \eta>0,\\
r\beta(0,C)C,&\qquad \textrm{ if } \eta=0. \end{array} \right.
\end{displaymath}
\end{lemma}

\begin{proof}
We assume that
\begin{displaymath}
\lim_{t \to +\infty} x^{\phi}(t)=C.
\end{displaymath}
By using (\ref{yexplicite}), we obtain that
\begin{displaymath}
y^{\psi}(t)=\int_{0}^r e^{-\eta
s}\beta(0,x^{\phi}(t-s))x^{\phi}(t-s)ds, \qquad \textrm{for }
t\geq r.
\end{displaymath}
Then, we easily conclude that
\begin{displaymath}
\lim_{t \to +\infty} y^{\psi}(t)=\bigg(\int_{0}^r e^{-\eta s} ds
\bigg)\beta(0,C)C.
\end{displaymath}
This ends the proof.
\end{proof}

Lemma \ref{lemmastability} implies that, in order to study the
stability of the solutions of Problem
(\ref{equationmatzero})-(\ref{equationmatzeroprime}), we only need
to concentrate on the stability of the solutions of the delay
differential equation (\ref{equationmatzero}).

In \cite{mackey1978}, Mackey has proposed that the function
$\beta(0,\cdot)$ is a Hill function, defined by
\begin{equation} \label{betahillfunction}
\beta(0,x)=\beta_{0}\frac{\theta^{n}}{\theta^{n}+x^{n}},
\end{equation}
where $\beta_{0}$ and $\theta$ are two positive constants and
$n\geq 1$. This function is used to describe, from a reasonable
biological point of view, the fact that the rate of re-entry in
the proliferating compartment is a decreasing function of the
total number of resting cells.

\noindent We recall that the function $x\mapsto\beta(0,x)$ is
supposed to be continuous and positive. From now on, we also
suppose that $x\mapsto\beta(0,x)$ is decreasing on $[0,+\infty)$
and satisfies
\begin{equation}\label{condlimbeta}
\lim_{x\to +\infty}\beta(0,x)=0.
\end{equation}
These assumptions have been done, for the first time, by Mackey
\cite{mackey1978} in 1978 and have been used by Mackey and
Rudnicki \cite{mackey1994} in 1994.

\noindent Before studying the stability of Problem
(\ref{equationmatzero})-(\ref{equationmatzeroprime}), we recall a
non-trivial property of the solutions of (\ref{equationmatzero}).
The result in Proposition \ref{propcondbounded} has been proved
for a similar equation by Mackey and Rudnicki \cite{mackey1994},
in 1994.

\begin{proposition} \label{propcondbounded}
Assume that $\rho >0$. Then, every solution of Equation
(\ref{equationmatzero}) is bounded.
\end{proposition}

One can notice that, if $\rho=0$, then Equation
(\ref{equationmatzero}) may have unbounded solutions. A counter
example is given in the next proposition.

\begin{proposition} \label{counterexample}
Assume that $\rho=0$ and that there exists $\overline{x}>0$ such
that the function $x\mapsto x\beta(0,x)$ is decreasing on
$[\overline{x}, +\infty)$. Let $\overline{\mu}\in C[0,1]$ and
$\Gamma\in C(\Omega)$ be such that (\ref{condcompatibilite})
holds, $\overline{\mu}(0)>\overline{x}$ and
\begin{equation}\label{condGamma}
2\xi(t,0)\Gamma(0,r-t) > \Gamma(0,0), \qquad \textrm{ for } t\in
[0,r].
\end{equation}
Then, the solution of Equation (\ref{equationmatzero}) is
unbounded.
\end{proposition}

\begin{proof}
Let consider the solution $x(t)$ of the problem
\begin{displaymath}
\left\{ \begin{array}{rcll}
x^{\prime}(t)&=& 2\overline{\xi}(0)\beta(0,x(t-r))x(t-r)-\beta(0,x(t))x(t),& \textrm{ for } t\geq r,\\
x(t)&=&\phi(t),& \textrm{ for } 0\leq t\leq r.
\end{array} \right.
\end{displaymath}
First, one can notice that, if $\lim_{t\to +\infty}x(t)$ exists
and is equal to $C$, then $C=0$.

\noindent By contradiction, if we suppose that $C>0$, then we
obtain that
\begin{displaymath}
\lim_{t\to +\infty}x^{\prime}(t) = \big(2\overline{\xi}(0) -1
\big)\beta(0,C)C>0,
\end{displaymath}
because Condition (\ref{condGamma}) implies that
$2\overline{\xi}(0)>1$. This contradicts the fact that $x(t)$
converges. Then $C=0$.

\noindent Secondly, let $\overline{\mu}\in C[0,1]$ and $\Gamma\in
C(\Omega)$ be such that (\ref{condcompatibilite}) and
(\ref{condGamma}) hold, and $\overline{\mu}(0)>\overline{x}$. The
solution $\phi(t)$ of the problem
\begin{displaymath}
\left\{ \begin{array}{rcll}
\phi^{\prime}(t)&=& 2\xi(t,0)\Gamma(0,r-t)-\beta(0,\phi(t))\phi(t),& \textrm{ for } t\geq r,\\
\phi(0)&=&\overline{\mu}(0),&
\end{array} \right.
\end{displaymath}
satisfies
\begin{displaymath}
\phi^{\prime}(0)= 2\xi(0,0)\Gamma(0,r)-\Gamma(0,0)>0.
\end{displaymath}
Consequently, there exists $\epsilon\in(0,r]$ such that
$\phi^{\prime}(t)>0$, for $t\in[0,\epsilon)$.

\noindent Hence, $\phi(0)<\phi(\epsilon)$. Then,
\begin{displaymath}
\begin{array}{rcl}
\phi^{\prime}(\epsilon)&>&\Gamma(0,0)-\beta(0,\phi(\epsilon))\phi(\epsilon),\\
 &\geq&\beta(0,\overline{\mu}(0))\overline{\mu}(0)-\beta(0,\phi(\epsilon))\phi(\epsilon),\\
 &\geq&0.
\end{array}
\end{displaymath}
By steps, we conclude that $\phi^{\prime}(t)>0$, for $t\in[0,r]$.

\noindent By the same way, we obtain that
\begin{displaymath}
x^{\prime}(r)>\beta(0,\phi(0))\phi(0)-\beta(0,\phi(r))\phi(r)\geq0.
\end{displaymath}
By using the same reasonning, we show that $x^{\prime}(t)>0$, for
$t\geq r$. Hence, $x$ is unbounded and the proof is complete.
\end{proof}

\begin{remark}
Even if the trivial solution of (\ref{equationmatzero}) is
unstable, the trivial solution of (\ref{equationmatzeroprime}) may
be stable. For example, if $\lim_{t\to +\infty}
x^{\phi}(t)=+\infty$, then, by using (\ref{yexplicite}) and
(\ref{condlimbeta}), we easily obtain that $\lim_{t\to +\infty}
y^{\psi}(t)=0$.
\end{remark}

\noindent The assumption on the function $x\mapsto x\beta(0,x)$,
in Proposition \ref{counterexample}, holds, for example, if
$\beta$ is given by (\ref{betahillfunction}) with $n>1$. In this
case, the function $x\mapsto x\beta(0,x)$ is decreasing for $x\geq
\overline{x}=\theta/(n-1)^{1/n}$.

\noindent We determine, in the next theorem, the global stability
area of the trivial solution of Equation (\ref{equationmatzero}).

\begin{theorem} \label{theoremstabilityzero}
The trivial solution of Equation (\ref{equationmatzero}) is
globally stable if and only if
\begin{displaymath}
(2\overline{\xi}(0)-1)\beta(0,0) < \rho.
\end{displaymath}
\end{theorem}

\begin{proof}
First, we assume that $(2\overline{\xi}(0)-1)\beta(0,0) <\rho$. We
are going to show the global stability by using a Lyapunov
functional.

\noindent We denote by $C^+$ the subset of $C[0,r]$ containing
nonnegative functions. We set
\begin{displaymath}
f(x)=x\beta(0,x)
\end{displaymath}
and
\begin{displaymath}
\mathcal{F}(x)=\displaystyle\int_{0}^{x} f(s)ds, \qquad \textrm{
for } x\geq0.
\end{displaymath}
We define the mapping $J:C^{+} \to \mathbb{R}$ by
\begin{displaymath}
J(\phi)=\mathcal{F}(\phi(r))+\overline{\xi}(0)\int_{0}^{r}f^{2}(\phi(\sigma))d\sigma,
\qquad \textrm{ for } \phi\in C^{+}.
\end{displaymath}
Then
\begin{displaymath}
\overset{\bullet}{J}(\phi) =
\frac{d\phi}{dt}(r)f(\phi(r))+\overline{\xi}(0)\big(f^{2}(\phi(r))-f^{2}(\phi(0))\big).\nonumber
\end{displaymath}
Since
\begin{displaymath}
\frac{d\phi}{dt}(r) =
-\rho\phi(r)-f\big(\phi(r)\big)+2\overline{\xi}(0)f\big(\phi(0)\big),
\end{displaymath}
then
\begin{displaymath}
\begin{array}{rcl}
\overset{\bullet}{J}(\phi) & = & -\rho\beta(0,\phi(r))\phi^2(r)-f^2(\phi(r))+\overline{\xi}(0)\Big(f^{2}(\phi(r))+ 2f(\phi(r))f(\phi(0))- f^{2}(\phi(0))\Big), \\
        & = & -\Big(\rho+\beta(0,\phi(r)) \Big)\beta(0,\phi(r))\phi^2(r) +2\overline{\xi}(0)f^{2}(\phi(r))-\overline{\xi}(0)\Big(f(\phi(r))- f(\phi(0))\Big)^2.
\end{array}
\end{displaymath}
Hence,
\begin{displaymath}
\overset{\bullet}{J}(\phi) \leq
-\Big(\rho-\big(2\overline{\xi}(0)-1\big)\beta(0,\phi(r))
\Big)\beta(0,\phi(r))\phi^2(r).
\end{displaymath}
Since $(2\overline{\xi}(0)-1)\beta(0,0) < \rho$ and the function
$x\mapsto \beta(0,x)$ is decreasing and positive on
$\mathbb{R}^+$, then the function
\begin{displaymath}
\lambda(u)=\Big(\rho-\big(2\overline{\xi}(0)-1\big)\beta(0,u)\Big)\beta(0,u)u^{2}
\end{displaymath}
is nonnegative on $\mathbb{R}^+$ and $\lambda(u)=0$ if and only if
$u=0$. Consequently, every solution of Equation
(\ref{equationmatzero}) with $\phi\in C^{+}$ tends to zero as $t$
tends to $+\infty$.

\noindent Now, if we assume that $\rho \leq
(2\overline{\xi}(0)-1)\beta(0,0)$, then, immediatly,
\begin{displaymath}
-(\rho+\beta(0,0))\geq -2\overline{\xi}(0)\beta(0,0).
\end{displaymath}
Hence, by using Bellman and Cooke (\cite{bellmancooke}, Theorem
13.8), we obtain that the trivial solution of
(\ref{equationmatzero}) is unstable.
\end{proof}

Remark that, by using Lemma \ref{lemmastability}, if
$(2\overline{\xi}(0)-1)\beta(0,0) <\rho$, then the trivial
solution of (\ref{equationmatzeroprime}) is also globally stable.

We are going to use the results of Theorem
\ref{theoremstabilityzero} in the next sections to obtain global
stability and instability for the solutions of Problem
(\ref{integratedformulationN})-(\ref{integratedformulationP}).

\section{Global stability for the maturity structured model} \label{asymptoticbehavioursection}

In this section, we establish a result of global stability for
Problem
(\ref{integratedformulationN})-(\ref{integratedformulationP})
which stresses the influence of immature cells on the total
population. First, we recall some definitions.

\begin{definition}
The trivial solution of Problem
(\ref{integratedformulationN})-(\ref{integratedformulationP}) is
\emph{locally stable} if, for all $\varepsilon>0$, there exist
$\nu>0$ and $T>0$ such that, if $\overline{\mu}\in C[0,1]$ and
$\Gamma\in C(\Omega)$ satisfy
\begin{displaymath}
\| \overline{\mu}\|<\nu \quad \textrm{ and } \quad \| \Gamma
\|_{\Omega}<\nu,
\end{displaymath}
then
\begin{equation}\label{defstabtrivialsolution}
\| N^{\overline{\mu},\Gamma}(t,.) \|<\varepsilon \quad \textrm{
and } \quad \| P^{\overline{\mu},\Gamma}(t,.) \|<\varepsilon,
\quad \textrm{ for } t\geq T.
\end{equation}
The trivial solution of Problem
(\ref{integratedformulationN})-(\ref{integratedformulationP}) is
\emph{globally stable} if, for all $\varepsilon>0$, there exists
$T>0$ such that (\ref{defstabtrivialsolution}) holds.
\end{definition}

Throughout this section, we assume that the function $x\mapsto
\beta(m,x)$ is uniformly bounded and that the function $x\mapsto
x\beta(m,x)$ is locally Lipschitz continuous for all $m\in[0,1]$.
In the next theorem, we show our main result, which makes the link
between the global stability of the trivial solution of Problem
(\ref{integratedformulationN})-(\ref{integratedformulationP}) and
the stability of the immature cells population.

\begin{theorem} \label{globalstability}
Assume that Condition (\ref{inegtauvarianlin}) holds. Let us
suppose that the trivial solution of Problem
(\ref{integratedformulationN})-(\ref{integratedformulationP}) is
locally stable. Then this solution is globally stable on the set
\begin{displaymath}
\Omega_{GS}=\Big\{ (\overline{\mu},\Gamma)\in C[0,1]\times
C(\Omega)\ \ : \ \ \lim_{t\to +\infty}
N^{\overline{\mu},\Gamma}(t,0)=\lim_{t\to +\infty}
P^{\overline{\mu},\Gamma}(t,0)=0 \Big\}.
\end{displaymath}
\end{theorem}

\begin{proof}
We first show that, if the trivial solution of Equation
(\ref{integratedformulationN}) is locally stable, then it is
globally stable on the set
\begin{displaymath}
\Omega_N=\Big\{ (\overline{\mu},\Gamma)\in C[0,1]\times C(\Omega)\
\ : \ \ \lim_{t\to +\infty} N^{\overline{\mu},\Gamma}(t,0)=0
\Big\}.
\end{displaymath}
Let us suppose that the trivial solution of
(\ref{integratedformulationN}) is locally stable. Then, for all
$\overline{\mu}\in C[0,1]$, $\Gamma\in C(\Omega)$ and
$\varepsilon>0$, there exist $\nu>0$ and $T>0$, such that, if
\begin{displaymath}
\| \overline{\mu}\|<\nu \quad \textrm{ and } \quad \| \Gamma
\|_{\Omega}<\nu,
\end{displaymath}
then,
\begin{equation} \label{condlocalstability}
|N^{\overline{\mu},\Gamma}(t,m)|<\varepsilon, \quad \textrm{ for }
t\geq T \textrm{ and } m\in[0,1].
\end{equation}
Let $\varepsilon>0$ be given and let
$(\overline{\mu},\Gamma)\in\Omega_N$. Then $\lim_{t\to +\infty}
N^{\overline{\mu},\Gamma}(t,0)=0$, so there exists $t_0>0$ such
that
\begin{displaymath}
|N^{\overline{\mu},\Gamma}(t,0)|<\frac{\nu}{2}, \quad \textrm{ for
} t\geq t_0.
\end{displaymath}
Let $\zeta\in C[0,1]$ and $\Upsilon\in C(\Omega)$ be given. Since
the solutions of Equation (\ref{integratedformulationN}) are
continuous, then there exists $\delta>0$ such that, if
\begin{displaymath}
|\zeta(m)-\overline{\mu}(0)|<\delta \quad \textrm{ and } \quad
|\Upsilon(m,a)-\Gamma(0,a)|<\delta,
\end{displaymath}
for $m\in[0,1]$ and $a\in[0,\tau_{max}]$, then
\begin{displaymath}
|N^{\zeta,\Upsilon}(t,m)-N^{\overline{\mu},\Gamma}(t,0)|<\frac{\nu}{2},
\end{displaymath}
for $t\in[t_0,t_0+\tau_{max}]$ and $m\in[0,1]$.

\noindent Now, since $\overline{\mu}$ and $\Gamma$ are continuous,
then there exists $b\in(0,1)$ such that
\begin{displaymath}
|\overline{\mu}(m)-\overline{\mu}(0)|<\delta \quad \textrm{ and }
\quad |\Gamma(m,a)-\Gamma(0,a)|<\delta,
\end{displaymath}
for $m\in[0,b]$ and $a\in[0,\tau_{max}]$. We define the following
functions,
\begin{displaymath}
\overline{\mu}_b(m)=\left\{ \begin{array}{ll}
\overline{\mu}(m),& \textrm{ if } m\in[0,b],\\
\overline{\mu}(b),& \textrm{ if } m\in[b,1],
\end{array}\right.
\quad \textrm{ and } \quad \Gamma_b(m,.)=\left\{ \begin{array}{ll}
\Gamma(m,.),& \textrm{ if } m\in[0,b],\\
\Gamma(b,.),& \textrm{ if } m\in[b,1].
\end{array}\right.
\end{displaymath}
Then, for $m\in[0,1]$ and $a\in[0,\tau_{max}]$, we get
\begin{displaymath}
|\overline{\mu}_b(m)-\overline{\mu}(0)|<\delta \quad \textrm{ and
} \quad |\Gamma_b(m,a)-\Gamma(0,a)|<\delta.
\end{displaymath}
Consequently,
\begin{displaymath}
|N^{\overline{\mu}_b,\Gamma_b}(t,m)-N^{\overline{\mu},\Gamma}(t,0)|<\frac{\nu}{2},
\end{displaymath}
for $t\in[t_0,t_0+\tau_{max}]$ and $m\in[0,1]$. It follows that
\begin{displaymath}
|N^{\overline{\mu}_b,\Gamma_b}(t,m)|<\nu, \quad \textrm{ for }
t\in[t_0,t_0+\tau_{max}] \textrm{ and } m\in[0,1].
\end{displaymath}

\noindent Since $N^{\overline{\mu}_b,\Gamma_b}(t,m)$ is a solution
of Equation (\ref{integratedformulationN}) and since this equation
is autonomous for $t$ large enough, then
$N^{\overline{\mu}_b,\Gamma_b}(t,m)$ becomes an initial condition
of (\ref{integratedformulationN}) on
$[t_0,t_0+\tau_{max}]\times[0,1]$. We deduce, from
(\ref{condlocalstability}), that there exists $\widetilde{T}\geq
t_0+\tau_{max}$ such that
\begin{displaymath}
|N^{\overline{\mu}_b,\Gamma_b}(t,m)|<\varepsilon, \qquad \textrm{
for } t\geq \widetilde{T} \textrm{ and } m\in[0,1].
\end{displaymath}

\noindent From Theorem \ref{smallmaturitycells}, there exists
$\overline{t}>0$ such that
\begin{displaymath}
N^{\overline{\mu}_b,\Gamma_b}(t,m)=N^{\overline{\mu},\Gamma}(t,m),
\qquad \textrm{ for } t\geq
\overline{t}+\tau_{max}-\ln\big(h(g(1))\big) \textrm{ and }
m\in[0,1].
\end{displaymath}
Hence,
\begin{displaymath}
 \|N^{\overline{\mu},\Gamma}(t,.)\|<\varepsilon, \qquad \textrm{ for } t\geq \max\big\{\widetilde{T},\overline{t}+\tau_{max}-\ln\big(h(g(1))\big)\big\}.
\end{displaymath}
Then, the trivial solution of Equation
(\ref{integratedformulationN}) is globally stable.

\noindent By the same way and using Corollary
\ref{smallmaturityprol}, we show that, if $\overline{\mu}$ and
$\Gamma$ are such that
\begin{displaymath}
\lim_{t\to +\infty} P^{\overline{\mu},\Gamma}(t,0)=0,
\end{displaymath}
then, the trivial solution of Equation
(\ref{integratedformulationP}) is globally stable. This completes
the proof.
\end{proof}

\begin{remark}\label{remark}
One has to notice that the result in Theorem \ref{globalstability}
allows us to obtain the global exponential stability of the
trivial solution of Problem
(\ref{integratedformulationN})-(\ref{integratedformulationP}) on
the set $\Omega_{GS}$, when this solution is locally exponentially
stable. The proof, in this case, is identical to the previous one.
\end{remark}

The behaviour of the immature cells population has been studied in
Section \ref{sectionbehaviourimmaturecells}. The local stability
of System
(\ref{integratedformulationN})-(\ref{integratedformulationP}) has
been studied by Adimy and Crauste \cite{adimycrauste}. The author
proved that, under the assumptions that the function $x\mapsto
\beta(m,x)$ is uniformly bounded and the function $x\mapsto
x\beta(m,x)$ is locally Lipschitz continuous for all $m\in[0,1]$,
then the trivial solution of
(\ref{integratedformulationN})-(\ref{integratedformulationP}) is
locally exponentially stable if
\begin{equation}\label{condlocstab}
\bigg(1+2\sup_{(m,t)\in\Omega_{\Delta}}\xi(t,m)\bigg)\sup_{m\in[0,1]}\Big(\beta(m,0)\Big)<\min\bigg\{
\inf_{m\in[0,1]} \Big(\delta(m)+V^{\prime}(m)\Big),
\inf_{m\in[0,1]} \Big(\gamma(m)+V^{\prime}(m)\Big) \bigg\}.
\end{equation}
The proof is based on an induction reasonning.

\noindent Then, we can deduce the following corollary, which deals
with the global stability of the system.

\begin{corollary}
Assume that Condition (\ref{inegtauvarianlin}) and Inequality
(\ref{condlocstab}) hold. Then the trivial solution of System
(\ref{integratedformulationN})-(\ref{integratedformulationP}) is
globally exponentially stable.
\end{corollary}

\begin{proof}
From Inequality (\ref{condlocstab}), we obtain the local
exponential stability of the trivial solution of
(\ref{integratedformulationN})-(\ref{integratedformulationP}).

\noindent From the definition of $\rho$, we get
\begin{displaymath}
\inf_{m\in[0,1]} \Big(\delta(m)+V^{\prime}(m)\Big)\leq \rho.
\end{displaymath}
Moreover,
\begin{displaymath}
\overline{\xi}(0) \leq \sup_{(m,t)\in\Omega_{\Delta}}\xi(t,m)
\quad \textrm{ and } \quad \beta(0,0)\leq
\sup_{m\in[0,1]}\beta(m,0).
\end{displaymath}
Hence, we obtain
\begin{displaymath}
(2\overline{\xi}(0)-1)\beta(0,0) < \rho.
\end{displaymath}
Then, Theorem \ref{theoremstabilityzero} yields to the global
stability of the trivial solution of
(\ref{equationmatzero})-(\ref{equationmatzeroprime})

\noindent By using Theorem \ref{globalstability} and Remark
\ref{remark}, we conclude.
\end{proof}

As an example, let us suppose that $\delta\geq0$ and $\gamma\geq0$
are constant, $V$ and $g$ are linear functions of the maturity
$m$, given, for $m\in[0,1]$, by
\begin{displaymath}
V(m)=m, \qquad \textrm{ and } \qquad g(m)=\frac{1}{\kappa}m, \quad
\textrm{ with } \kappa>1,
\end{displaymath}
and the function $\beta$ is a Hill function (see Mackey
\cite{mackey1978}), defined by
\begin{displaymath}
\beta(m,x)=\beta_0(m)\frac{\theta^n(m)}{\theta^n(m)+x^n},
\end{displaymath}
with $\beta_0$ and $\theta$ two continuous and positive functions
on $[0,1]$, and $n>1$. One can remark that, in this case, the
function $V$ satisfies Condition (\ref{(H.1).1}).

\noindent Furthermore, we assume that the function $\tau$ is
given, for $m\in[0,1]$, by
\begin{displaymath}
\tau(m)=\ln(m+\alpha), \qquad \textrm{ with } \alpha>1.
\end{displaymath}
In this case, $\tau$ is increasing. Therefore, Condition
(\ref{condtau}) is satisfied. We obtain that
\begin{displaymath}
\Delta(m)=\frac{1}{2}\Big(\sqrt{4\kappa m +\alpha^2} -\alpha
\Big), \qquad \textrm{ for } m\in[0,1],
\end{displaymath}
and the characteristic curves are given by
\begin{displaymath}
\chi(s,m)=me^s, \qquad \textrm{ for } s\leq0 \textrm{ and }
m\in[0,1].
\end{displaymath}
\noindent  By remarking that $\Delta(m)<m$ for $m\in(0,1]$ if and
only if $\alpha>\kappa$, then we obtain that the trivial solution
of (\ref{integratedformulationN})-(\ref{integratedformulationP})
is globally exponentially stable if
\begin{displaymath}
(1+2\kappa)\sup_{m\in[0,1]}\beta_0(m)<\min\{ \delta, \gamma\}
\qquad \textrm{ and } \qquad \alpha>\kappa.
\end{displaymath}

In the next section, we conclude our asymptotic study by giving a
result of instability, based on the results of Section
\ref{sectionbehaviourimmaturecells}.

\section{Instability}\label{sectioninstability}

The trivial solution of (\ref{integratedformulationN}) is
\emph{unstable} if it is not stable, this means, if there exists $
\varepsilon>0$ such that, for all $\nu>0$, there exists
$(\overline{\mu}, \Gamma)\in C[0,1]\times C(\Omega)$ which
satisfies
\begin{displaymath}
\| \overline{\mu}\|<\nu \quad \textrm{ and } \quad \| \Gamma
\|_{\Omega}<\nu,
\end{displaymath}
and
\begin{displaymath}
\| N^{\overline{\mu},\Gamma}(t,.) \|>\varepsilon, \qquad \textrm{
for } t\geq 0.
\end{displaymath}
In the next theorem, we show that the instability of the immature
cells population leads to the instability of the entire
population.

\begin{theorem}
Assume that
\begin{equation}\label{condinstability}
\rho \leq (2\overline{\xi}(0)-1)\beta(0,0).
\end{equation}
Then, the trivial solution of Problem
(\ref{integratedformulationN}) is unstable.
\end{theorem}

\begin{proof}
From Theorem \ref{theoremstabilityzero} and
(\ref{condinstability}), we obtain that the trivial solution of
Equation (\ref{equationmatzero}) is unstable. That is there exist
$\varepsilon>0$, $\overline{\mu}\in C[0,1]$ and $\Gamma\in
C(\Omega)$ such that $N^{\overline{\mu},\Gamma}(t,0)$ does not
tend to zero when $t$ goes to infinity. Then, there exist
$\varepsilon>0$ and $(t_n)_{n\in\mathbb{N}}$, with $t_n \to
+\infty$, such that
\begin{displaymath}
N^{\overline{\mu},\Gamma}(t_n,0)>\varepsilon, \quad \textrm{ for }
n\in\mathbb{N}.
\end{displaymath}
Let us suppose, by contradiction, that the trivial solution of
(\ref{integratedformulationN}) is stable. Then, in particular,
there exist $\nu>0$ and $T>0$ such that, if
\begin{displaymath}
\| \overline{\mu}\|<\nu \quad \textrm{ and } \quad \| \Gamma
\|_{\Omega}<\nu,
\end{displaymath}
then
\begin{displaymath}
\| N^{\overline{\mu},\Gamma}(t,.) \|<\varepsilon, \quad \textrm{
for } t\geq T.
\end{displaymath}
Consequently,
\begin{displaymath}
| N^{\overline{\mu},\Gamma}(t_n,0) |<\varepsilon, \quad \textrm{
for } n\in\mathbb{N} \textrm{ such that } t_n\geq T.
\end{displaymath}
Since we can choose $\overline{\mu}$ and $\Gamma$ as small as
necessary, this yields a contradiction. We deduce the instability
of the trivial solution of (\ref{integratedformulationN}).
\end{proof}

One can remark that, even if the trivial solution of
(\ref{integratedformulationN}) is unstable, the trivial solution
of (\ref{integratedformulationP}) may be stable.


\begin{thebibliography}{99}

%
\bibitem{adimycrauste2003} M. Adimy and F. Crauste, \emph{Global stability of a partial differential equation with distributed delay due to cellular replication}, Nonlinear Analysis {\bf 54}, 1469-1491 (2003).
%
\bibitem{adimycrauste} M. Adimy and F. Crauste, \emph{Existence, positivity and stability for a model of cellular proliferation}, submitted.
%
\bibitem{adimypujo} M. Adimy and L. Pujo-Menjouet, \emph{A singular transport model describing cellular division}, C. R. Acad. Sci., Paris, Ser. I, Math {\bf 332}, 12, 1071-1076 (2001).
%
\bibitem{adimypujo2} M. Adimy and L. Pujo-Menjouet, \emph{Asymptotic behaviour of a singular transport equation modelling cell division}, Dis. Cont. Dyn. Sys. Ser. B \textbf{3}, 439-456 (2003).
%
\bibitem{adimypujo2003} M. Adimy and L. Pujo-Menjouet, \emph{A mathematical model describing cellular division with a proliferating phase duration depending on the maturity of cells}, submitted.
%
\bibitem{bellmancooke} R. Bellman and K.L. Cooke, \emph{Differential difference equations}, Academic Press, New-York-London, 462 p., 1963.
%
\bibitem{bradford} G. Bradford, B. Williams, R. Rossi and I. Bertoncello, \emph{Quiescence, cycling, and turnover in the primitive haematopoietic stem cell compartment}, Exper. Hematol. {\bf 25}, 445-453 (1997).
%
\bibitem{burns} F.J. Burns and I.F. Tannock, \emph{On the existence of a $G_{0}$ phase in the cell cycle}, Cell. Tissue Kinet. {\bf 19}, 321-334 (1970).
%
\bibitem{crabb1996_1} R. Crabb, J. Losson and M.C. Mackey, \emph{Dependence on initial conditions in non local PDE's and heredetary dynamical systems}, Proc. Inter. Conf. Nonlin. Anal. {\bf 4} (Tampa Bay, de Gruyter, Berlin, 1996) 3125-3136.
%
\bibitem{crabb1996_2} R. Crabb, M.C. Mackey and A. Rey, \emph{Propagating fronts, chaos and multistability in a cell replication model}, Chaos {\bf 6}, 477-492 (1996).
%
\bibitem{webb1996} J. Dyson, R. Villella-Bressan and G.F. Webb, \emph{A singular transport equation modelling a proliferating maturity structured cell population}, Can. Appl. Math. Quart. {\bf 4}, 65-95 (1996).
%
\bibitem{webb2000} J. Dyson, R. Villella-Bressan and G.F. Webb,\emph{ A nonlinear age and maturity structured model of population dynamics. I : Basic theory.}, J. Math. Anal. Appl. {\bf 242}, 1, 93-104 (2000).
%
\bibitem{webb2000_2} J. Dyson, R. Villella-Bressan and G.F. Webb, \emph{A nonlinear age and maturity structured model of population dynamics. II : Chaos.}, J. Math. Anal. Appl. {\bf 242}, 2, 255-270 (2000).
%
\bibitem{halelunel} J. Hale and S.M. Verduyn Lunel, {\em Introduction to functional differential equations\/}, Applied Mathematical Sciences 99, Springer-Verlag, New York, 447 p., 1993.
%
\bibitem{john} P.C.L John, \emph{The cell cycle}, London, Cambridge University Press, 1981.
%
\bibitem{mackey1978} M.C. Mackey, \emph{Unified hypothesis of the origin of aplastic anaemia and periodic hematopoiesis}, Blood {\bf 51},  941-956 (1978).
%
\bibitem{mackey1993} M.C. Mackey and A. Rey,  \emph{Multistability and boundary layer development in a transport equation with retarded arguments}, Can. Appl. Math. Quart. {\bf 1}, 1-21 (1993).
%
\bibitem{mackey1995_1} M.C. Mackey and A. Rey, \emph{Transitions and kinematics of reaction-convection fronts in a cell population model}, Physica D {\bf 80}, 120-139 (1995).
%
\bibitem{mackey1995_2} M.C. Mackey and A. Rey, \emph{Propagation of population pulses and fronts in a cell replication problem: non-locality and dependence on the initial function}, Physica D {\bf 86}, 373-395 (1995).
%
\bibitem{mackey1994} M.C. Mackey and R. Rudnicki, \emph{Global stability in a delayed partial differential equation describing cellular replication}, J. Math. Biol. {\bf 33}, 89-109 (1994).
%
\bibitem{mackey1999} M.C. Mackey and R. Rudnicki, \emph{A new criterion for the global stability of simultaneous cell replication and maturation processes}, J. Math. Biol. {\bf 38}, 195-219 (1999).
%
\bibitem{mitchison} J.M Mitchison, \emph{The biology of the cell cycle}, London, Cambridge University Press, 1971.
%
\bibitem{sachs} Sachs L., \emph{The molecular control of hemopoiesis and leukemia}, C. R. Acad. Sci. Paris \textbf{316}, 882-891 (1993).
%
\end{thebibliography}
\end{document}